\documentclass{article}
\usepackage{graphicx} 
\usepackage{amssymb}  
\usepackage{amsthm}
\usepackage{amsmath}
\usepackage{hyperref}

\title{Log-Linear Reaction Quotient Dynamics}
\author{Steven Diamond \\ steven@gridmatic.com \\ Gridmatic}

\begin{document}

\maketitle

\begin{abstract}
Chemical reaction networks in living cells maintain precise control over 
thousands of metabolites
despite operating far from equilibrium under constant perturbations.
While mass action kinetics accurately describe the underlying dynamics,
the resulting nonlinear differential equations are difficult to analyze and 
control,
particularly for large networks.
We propose a simplified model where reaction quotients (the ratios that measure how far reactions 
are from equilibrium)
evolve exponentially toward their equilibrium values when viewed on a 
logarithmic scale.
This principle leads to linear dynamics in log-space, providing several 
key advantages:
analytical solutions exist for arbitrary network topologies,
thermodynamic constraints are automatically satisfied through the relationship 
between
reaction quotients and Gibbs free energy,
conservation laws decouple from reaction quotient dynamics simplifying both 
analysis and control design,
and external energy sources couple linearly to the dynamics, unifying diverse 
biological regulatory mechanisms.
\end{abstract}

\section{Introduction}

Chemical reaction networks govern the behavior of living systems, from 
metabolic pathways to gene regulatory circuits 
\cite{Alon2006,Palsson2015,Heinrich1996}. Understanding and controlling these 
networks remains a central challenge in systems biology and synthetic biology 
\cite{DelVecchio2016,Khalil2010,Nielsen2016}. While mass action kinetics 
accurately describe the underlying dynamics 
\cite{Waage1864,Epstein1998,Erdi1989}, the resulting nonlinear differential 
equations are difficult to analyze, especially for large networks 
\cite{Gunawardena2003,Craciun2006}.

\paragraph{The control problem in biology.} 
Living cells maintain precise control over thousands of simultaneous chemical 
reactions, keeping metabolite concentrations far from equilibrium despite 
constant perturbations \cite{Fell1997,Sauro2004,Qian2007}. This regulation 
occurs through allosteric feedback \cite{Monod1965,Koshland1966}, covalent 
modification \cite{Goldbeter1990,Ferrell1996}, and transcriptional control 
\cite{Ptashne2002,Bintu2005}. How cells achieve robust control with such 
complex networks remains poorly understood \cite{Barkai1997,Yi2000,ElSamad2005}.

\paragraph{Existing frameworks.} 
To address the challenge of understanding metabolic regulation, several 
mathematical frameworks have been developed, each with distinct limitations. 
Chemical reaction network theory provides qualitative results about 
multistability and persistence \cite{FeinbergCRNT,HornJackson1972} but not 
quantitative dynamics. Metabolic control analysis quantifies sensitivity 
coefficients \cite{Kacser1973,Heinrich1974} but operates only near steady 
state. Flux balance analysis predicts metabolic fluxes 
\cite{Varma1994,Edwards2001,Orth2010} including with thermodynamic constraints 
\cite{Beard2002,Henry2007,Schellenberger2011}, but remains fundamentally 
static. Biochemical systems theory uses power-law approximations that become 
linear in logarithmic coordinates \cite{Savageau1969,Voit2000,Torres2002}, 
similar to our approach but using concentrations rather than reaction quotients 
as state variables.

\paragraph{Thermodynamic foundation.} 
The nonequilibrium thermodynamics of metabolism has long recognized that 
thermodynamic forces drive metabolic fluxes \cite{Onsager1931,Prigogine1955}. 
The relationship between reaction rates and Gibbs free energy has been explored 
extensively \cite{Noor2013,Park2016,Beard2008}. Our key innovation is using 
reaction quotients $Q$ as state variables, which naturally incorporates 
thermodynamic constraints through $\Delta G = RT\ln(Q/K_{eq})$ 
\cite{Alberty2003}. Building on this thermodynamic foundation, we develop a 
framework that makes reaction network dynamics analytically tractable

\paragraph{Our contribution.} 
We introduce a mathematical framework where reaction quotients follow 
first-order dynamics in logarithmic space: $d\ln Q/dt = -k\ln(Q/K_{eq})$. This 
yields: (1) analytical solutions for arbitrary reaction networks, (2) automatic 
incorporation of thermodynamic constraints, (3) decoupling of reaction dynamics 
from conservation laws, and (4) linear coupling to external energy sources. We 
demonstrate these advantages through examples including feedback inhibition, 
ATP-driven reactions, coupled transport, and glycolytic oscillations. The 
framework provides a tractable alternative to mass action kinetics while 
maintaining biological relevance.

The outline of our paper is as follows. Section 2 describes our framework in 
the context of a single chemical reaction. Section 3 extends our framework to 
systems of chemical reactions. Section 4 presents numerical examples. Section 5 
concludes.

\section{Single chemical reaction}

\paragraph{Dynamics.} For the reaction $aA + bB \rightleftharpoons cC + dD$, 
the reaction quotient $Q$ is defined as
\[
Q = \frac{[C]^c [D]^d}{[A]^a [B]^b}.
\]
We propose a mathematical approximation of chemical dynamics based on the 
principle that 
in log-space the reaction quotient relaxes toward equilibrium at a rate 
proportional to its distance from $K_{eq}$. This leads to the differential 
equation
\begin{equation}\label{eq:diff-eq}
\frac{d\ln Q}{dt} = -k \ln \left( Q/K_{eq}\right),
\end{equation}
where $k$ is the relaxation rate constant.

Taking the initial condition $Q_0 > 0$, the solution to the differential 
equation~\eqref{eq:diff-eq} is given by
\[
Q(t) = K_{eq}  \left(\frac{Q_0}{K_{eq}}\right)^{e^{-kt}},
\]
or equivalently in log-space
\[
\ln\left(\frac{Q(t)}{K_{eq}}\right) = \ln\left(\frac{Q_0}{K_{eq}}\right)e^{-kt}.
\]
The time constant $\tau = 1/k$ characterizes the exponential relaxation of 
$\ln(Q/K_{eq})$ towards zero, 
with the system reaching $(1-1/e) \approx 63\%$ of its equilibrium value on a 
logarithmic scale after time $\tau$.

\paragraph{Comparison with mass action.} The proposed reaction quotient 
dynamics differ from mass action.
Consider the reaction $A \rightleftharpoons B$ with rate constants $k_f$ and 
$k_r$. The reaction quotient dynamics derived from mass action are given by
\begin{equation}\label{eq:mass-action-q}
\frac{dQ}{dt} = k_r(1 + Q)(K_{eq} - Q).
\end{equation}
%
%
%
%
%
%
%
%
Near equilibrium, the mass action dynamics~\eqref{eq:mass-action-q} are well 
approximated by the linear dynamics
\[
\frac{dQ}{dt} \approx k_r (1 + K_{eq}) (K_{eq} - Q).
\]
The reaction quotient dynamics~\eqref{eq:diff-eq} are well approximated near 
equilibrium by
\[
\frac{dQ}{dt} \approx k (K_{eq} - Q).
\]
If we set $k =  k_r (1 + K_{eq})$, the two dynamics match near equilibrium. 

\paragraph{Conservation laws.} Conservation laws such as mass balance 
constraints are not enforced explicitly in the reaction quotient dynamics.
Instead, the dynamics leave absolute concentrations underdetermined. 
Conservation laws impose additional constraints that can be used to derive 
absolute concentrations.
In the reaction $A \rightleftharpoons B$ for example, the absolute 
concentrations $[A]$ and $[B]$ can be derived from $Q$ and the total 
concentration $C_{\text{total}}$:
\begin{align*}
[A] &= \frac{C_{\text{total}}}{1 + Q} \\
[B] &= \frac{C_{\text{total}}Q}{1 + Q}.
\end{align*}
The fact that reaction quotients are not constrained by conservation laws is 
well known in chemical reaction network theory 
\cite{FeinbergCRNT,HornJackson1972}.
We include a self-contained proof in Appendix~\ref{app:proof}.

\paragraph{Thermodynamic interpretation.} 
The reaction quotient dynamics~\eqref{eq:diff-eq} can be understood from a 
thermodynamic perspective. The reaction quotient $Q$ and equilibrium constant 
$K_{eq}$ are related to the Gibbs free energy through
\[
\Delta G = \Delta G^\circ + RT \ln Q = RT \ln \left(\frac{Q}{K_{eq}}\right),
\]
where $\Delta G^\circ = -RT \ln K_{eq}$ is the standard Gibbs free energy 
change. Thus, the reaction quotient dynamics~\eqref{eq:diff-eq} can be 
rewritten as
\begin{equation}\label{eq:gibbs-dynamics}
\frac{d(\Delta G)}{dt} = -k\Delta G.
\end{equation}
This shows that the Gibbs free energy decays exponentially toward equilibrium 
with time constant $\tau = 1/k$, consistent with the principle that spontaneous 
processes minimize free energy.

\paragraph{Control input.}
We hope not only to model a reaction's movement towards equilibrium but also 
maintain the reaction at non-equilibrium points. To maintain the reaction away 
from equilibrium we need a control input. We model the control input $u$ as 
linear in the reaction quotient dynamics:
\begin{equation}\label{eq:diff-eq-control-input}
\frac{d\ln Q}{dt} = -k \ln \left( Q/K_{eq}\right) + u.
\end{equation}
The solution to the differential equation~\eqref{eq:diff-eq-control-input} is 
given by
\[
Q(t) = K_{eq} \exp\left\{\left[\ln\left(\frac{Q_0}{K_{eq}}\right) - 
\frac{u}{k}\right]e^{-kt} + \frac{u}{k}\right\}.
\]
As $t \to \infty$, $Q$ converges to the steady state 
\[
Q_{ss} = K_{eq} e^{u/k}.
\]
From a thermodynamic perspective, the control input represents an external 
energy gradient coupled to the reaction:
\[
u = k_u  \Delta E/(RT),
\]
where $k_u$ is the coupling rate 
and $\Delta E$ is the external energy source.
At steady state, the system reaches a balance between the reaction's intrinsic 
thermodynamic driving force and the external energy input.

\section{Systems of chemical reactions}\label{sec:systems}

\paragraph{Reaction quotient vector dynamics.} Consider a system of $r$ coupled 
chemical reactions with reaction quotients $Q = (Q_1, Q_2, \ldots, Q_r) \in 
\mathbb{R}_{>0}^r$ and equilibrium constants $K_{eq} = (K_{eq,1}, K_{eq,2}, 
\ldots, K_{eq,r}) \in \mathbb{R}_{>0}^r$. We extend the single-reaction 
dynamics to the multivariable case:
\[
\frac{d\ln Q_i}{dt} = -\sum_{j=1}^r K_{ij} \ln(Q_j/K_{eq,j}), \quad i = 1, 
\ldots, r,
\]
where $K \in \mathbb{R}^{r \times r}$ is the relaxation rate matrix. In vector 
form:
\[
\frac{d}{dt}\ln Q = -K \ln(Q/K_{eq}),
\]
where the logarithm and division are applied element-wise.

\paragraph{Change of variables.} Define $x = \ln(Q/K_{eq}) \in \mathbb{R}^r$. 
The dynamics become linear:
\[
\frac{dx}{dt} = -Kx.
\]
For $K$ with positive eigenvalues, the solution is
\[
x(t) = e^{-Kt}x_0,
\]
which in terms of reaction quotients gives
\[
Q(t) = K_{eq} \circ \exp(e^{-Kt}\ln(Q_0/K_{eq})),
\]
where $\circ$ denotes element-wise multiplication and $\exp(\cdot)$ is applied 
element-wise.

\paragraph{Stability and equilibrium.} If all eigenvalues of $K$ have positive 
real parts, then $x(t) \to 0$ as $t \to \infty$, implying $Q(t) \to K_{eq}$. 
The rate of convergence is determined by the smallest eigenvalue of $K$.

\paragraph{Diagonal relaxation.} When $K = \text{diag}(k_1, k_2, \ldots, k_r)$ 
is diagonal, the reactions evolve independently:
\[
Q_i(t) = K_{eq,i} \left(\frac{Q_{0,i}}{K_{eq,i}}\right)^{e^{-k_i t}}, \quad i = 
1, \ldots, r.
\]
Each reaction relaxes to equilibrium with its own time constant $\tau_i = 
1/k_i$.

\paragraph{Coupled relaxation.} Off-diagonal elements in $K$ couple the 
reaction dynamics. For symmetric positive definite $K$, the system can be 
decomposed into eigenmodes:
\[
K = V\Lambda V^T,
\]
where $V$ contains the orthonormal eigenvectors and $\Lambda = 
\text{diag}(\lambda_1, \ldots, \lambda_r)$ contains the eigenvalues. In the 
eigenmode coordinates $z = V^T x$, the dynamics decouple:
\[
\frac{dz_i}{dt} = -\lambda_i z_i, \quad i = 1, \ldots, r.
\]

\paragraph{Control input.} We extend the controlled dynamics to multiple 
reactions:
\[
\frac{d\ln Q}{dt} = -K \ln(Q/K_{eq}) + u,
\]
where $u \in \mathbb{R}^r$ is the control input vector. In the transformed 
coordinates:
\[
\frac{dx}{dt} = -Kx + u.
\]
For constant $u$, the steady-state solution is
\[
x_{ss} = K^{-1}u,
\]
giving steady-state reaction quotients
\[
Q_{ss} = K_{eq} \circ \exp(K^{-1}u).
\]
Here $\circ$ denotes element-wise multiplication and $\exp(\cdot)$ is applied 
element-wise.

\paragraph{Thermodynamic interpretation.} Each component of the control vector 
represents an external energy gradient coupled to its respective reaction:
\[
u_i = k_{u,i} \Delta E_i/(RT),
\]
where $k_{u,i}$ is the coupling rate for reaction $i$ and $\Delta E_i$ is the 
external energy source. At steady state, the system balances the intrinsic 
thermodynamic forces with the external energy inputs across all reactions.

\paragraph{Conservation constraints.} As discussed in Appendix~\ref{app:proof}, 
conservation laws do not constrain the achievable reaction quotient states. For 
any $Q \in \mathbb{R}_{>0}^r$ satisfying $\ln Q \in \text{Im}(S^T)$ where $S$ 
is the stoichiometric matrix, there exist concentrations $c > 0$ that yield 
these reaction quotients while satisfying all mass balance constraints. 
Additionally, the constraint $\ln Q \in \text{Im}(S^T)$ automatically ensures 
that Wegscheider cycles are thermodynamically consistent: for any closed cycle 
of reactions, the product of equilibrium constants around the cycle equals one, 
preventing violations of detailed balance. This 
separation between reaction quotient dynamics and conservation laws simplifies 
both analysis and control design.





\section{Numerical examples}

We illustrate the log-linear framework through stylized examples that 
demonstrate key dynamical behaviors. While we focus on simple systems to build 
intuition, future work will validate the framework against experimental 
metabolic data. All code for reproducing these examples and exploring parameter 
variations is available at

\begin{center}
\url{https://github.com/reactionquotient/examples}
\end{center}

\begin{figure}[h]
\centering
\includegraphics[width=\textwidth]{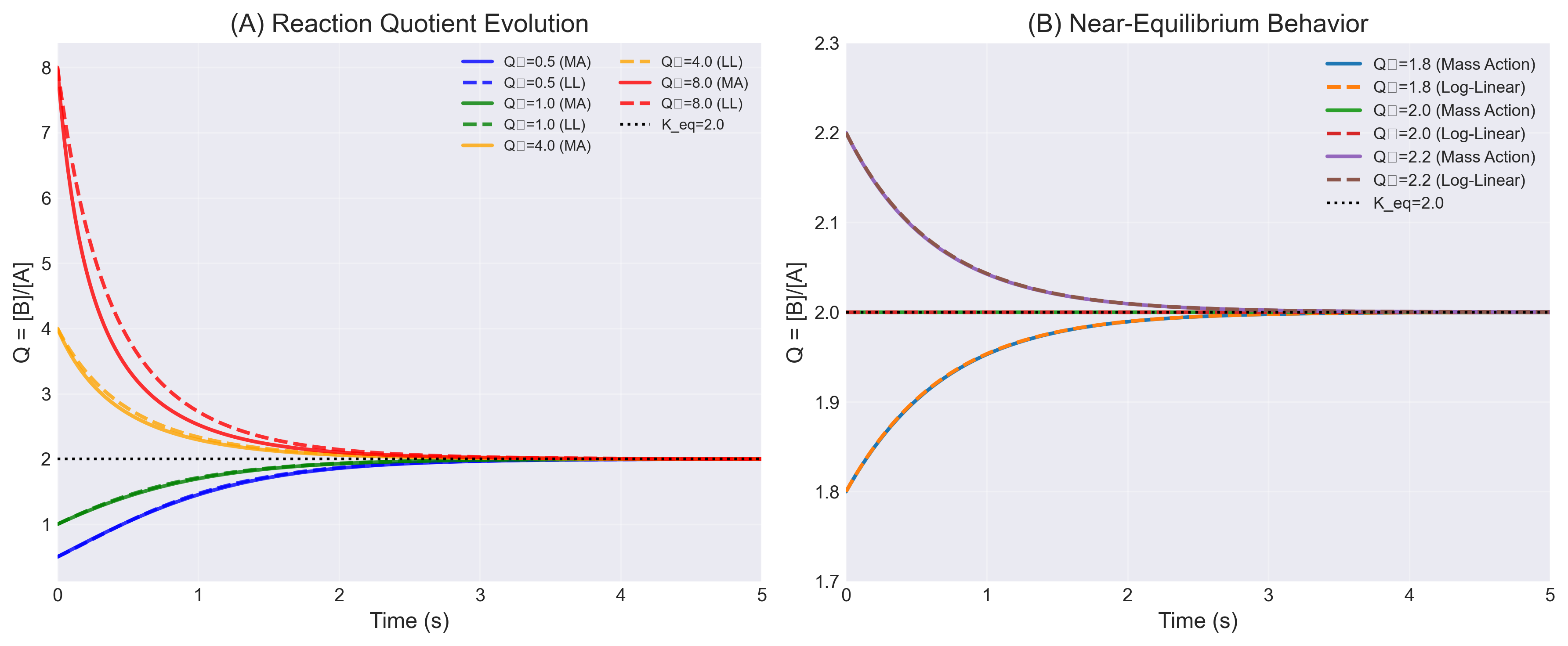}
\caption{Mass action (solid) vs log-linear (dashed) dynamics for $A 
\rightleftharpoons B$. (A) Evolution from various initial conditions. (B) 
Near-equilibrium behavior showing excellent agreement.}
\label{fig:AB-reaction}
\end{figure}

\paragraph{Mass action comparison.} We compare log-linear and mass action 
dynamics for $A \rightleftharpoons B$ with $k_f = 1.0$ s$^{-1}$ and $K_{eq} = 
2.0$. Setting $k = k_r(1 + K_{eq}) = 1.5$ s$^{-1}$ matches the dynamics near 
equilibrium. Figure~\ref{fig:AB-reaction} shows evolution from various initial 
conditions $Q_0 \in \{0.5, 1.0, 4.0, 8.0\}$. Both models converge to $Q = 
K_{eq}$, with similar trajectories despite different mathematical forms. Far 
from equilibrium (e.g., $Q_0 = 8$), the log-linear model relaxes faster due to 
exponential decay in log-space. Near equilibrium (panel B), the models are 
virtually indistinguishable, validating our parameter matching.

\begin{figure}[h]
\centering
\includegraphics[width=\textwidth]{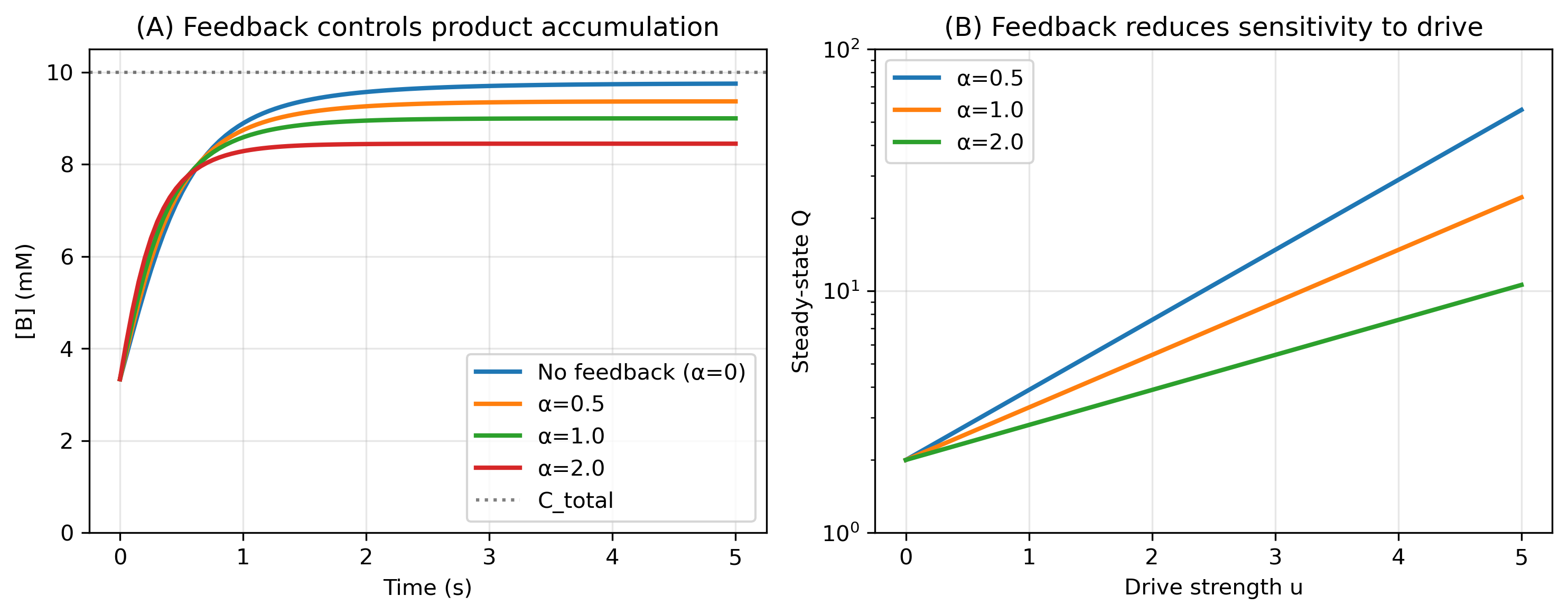}
\caption{Feedback inhibition in the log-linear framework. (A) Time evolution of 
$[B]$ under constant drive $u=3$ for different feedback strengths. Without 
feedback, $[B] \to C_{\text{total}}$; with feedback, accumulation is limited. 
(B) Steady-state response curves showing reduced sensitivity to drive as 
feedback increases.}
\label{fig:feedback}
\end{figure}

\paragraph{Feedback inhibition.}
Consider $A \to B$ where product $B$ inhibits its own formation, the minimal 
feedback motif. With $Q = [B]/[A]$ and conservation $[A] + [B] = 
C_{\text{total}}$, the dynamics become 
\[
\frac{d\ln Q}{dt} = -(k+\alpha)\ln(Q/K_{eq}) + u,
\]
where $\alpha$ is feedback strength and $u$ represents external drive. The 
feedback effectively increases the relaxation rate from $k$ to $k+\alpha$. The 
steady state 
\[
Q_{ss} = K_{eq}\exp(u/(k+\alpha))
\]
shows how feedback prevents unlimited accumulation. Figure~\ref{fig:feedback}A 
demonstrates this effect with $u=3$: without feedback ($\alpha=0$), $[B]$ 
approaches $C_{\text{total}}$ as all $A$ converts to $B$, while increasing 
feedback strength progressively limits $[B]$ accumulation despite constant 
drive. Figure~\ref{fig:feedback}B shows the input-output relationship 
$Q_{ss}(u)$ for different $\alpha$ values, revealing how stronger feedback 
reduces sensitivity to drive fluctuations. The system remains stable for all 
$\alpha > 0$, implementing a natural robustness-responsiveness tradeoff.

\begin{figure}[h]
\centering
\includegraphics[width=\textwidth]{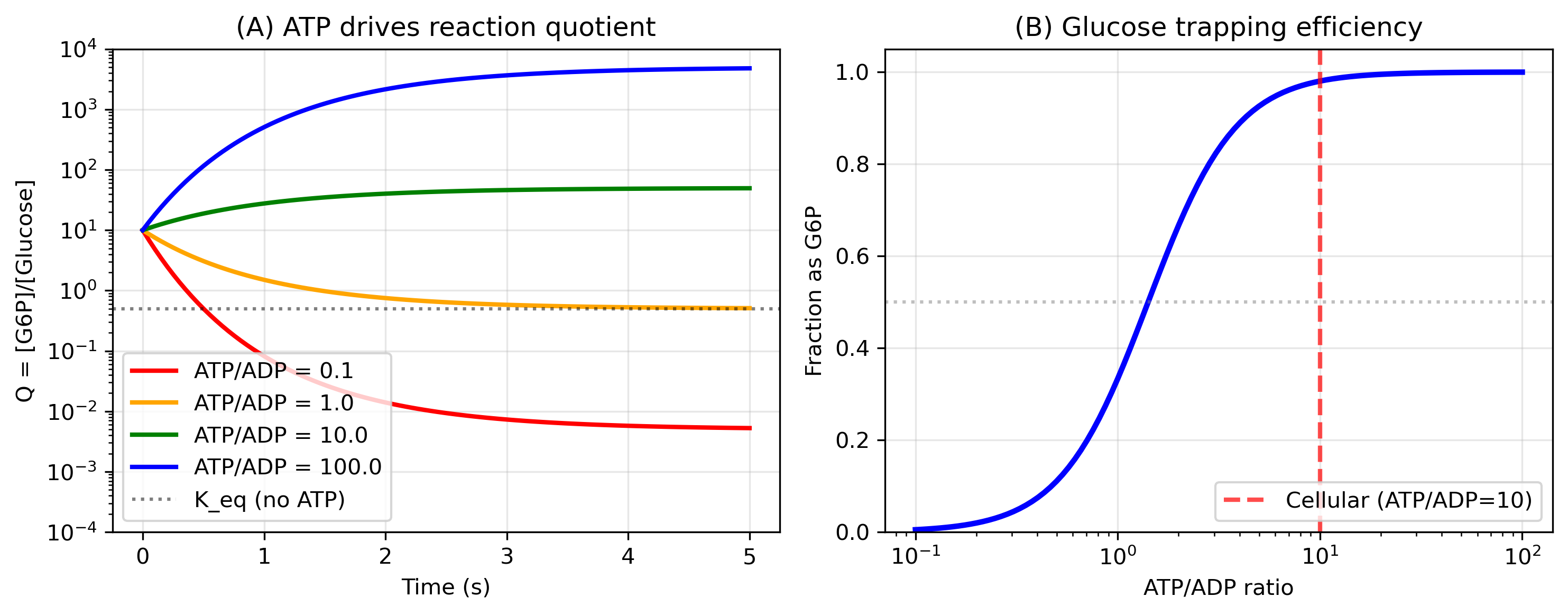}
\caption{ATP-driven glucose phosphorylation by hexokinase. (A) Reaction 
quotient dynamics for different ATP/ADP ratios. Low energy (red) reverses the 
reaction; high energy (blue) drives strong phosphorylation. (B) Sigmoid 
relationship between cellular energy state and glucose trapping efficiency. The 
cellular ratio (red dashed line) operates in the sensitive region.}
\label{fig:hexokinase}
\end{figure}

\paragraph{ATP-driven reaction.}
Consider the hexokinase reaction, the first step of glycolysis where glucose is 
phosphorylated using ATP. We model this as $Q = [\text{G6P}]/[\text{Glucose}]$ 
with the cellular ATP/ADP ratio providing external drive. With $K_{eq} = 0.5$ 
(chemistry alone favors glucose), conservation $[\text{Glucose}] + [\text{G6P}] 
= C_{\text{total}}$, and ATP coupling strength $k_{\text{ATP}} = 2$, the 
dynamics become
\[
\frac{d\ln Q}{dt} = -k\ln(Q/K_{eq}) + 
k_{\text{ATP}}\ln([\text{ATP}]/[\text{ADP}]).
\]
Figure~\ref{fig:hexokinase}A shows how different metabolic states affect the 
reaction: low ATP/ADP = 0.1 (starvation) drives the reaction backward to 
release glucose, while normal cellular ATP/ADP = 10 drives $Q$ to 50, two 
orders of magnitude above chemical equilibrium. Figure~\ref{fig:hexokinase}B 
reveals a sharp sigmoid response—the glucose trapping efficiency jumps from 
near zero to nearly complete between ATP/ADP ratios of 1 to 10. At cellular 
conditions, 98\% of glucose is trapped as G6P despite chemistry favoring the 
reverse. This example shows that log-linear reaction quotient dynamics can 
capture biological switching behavior: ATP-dependent control emerges directly 
from the thermodynamic coupling term, producing the sharp sigmoid response that 
cells exploit for metabolic regulation.

\begin{figure}[h]
\centering
\includegraphics[width=\textwidth]{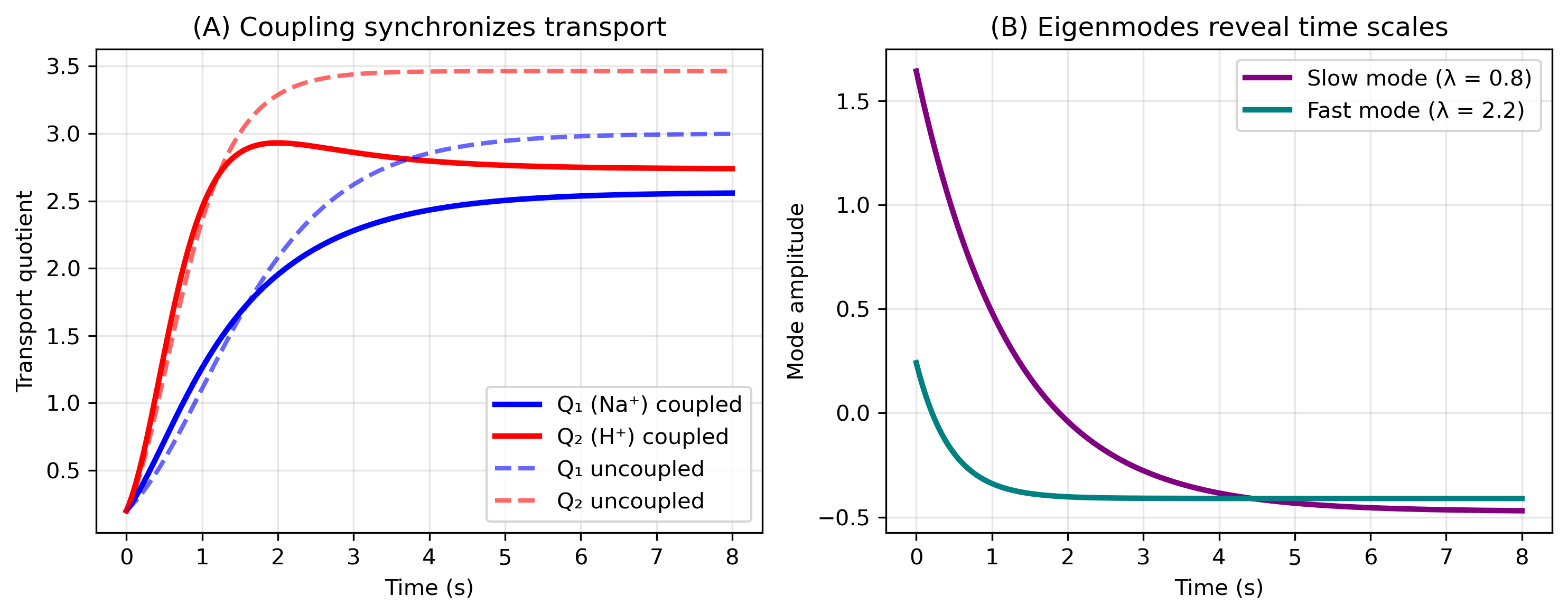}
\caption{Coupled membrane transport. (A) Off-diagonal coupling creates 
overshoot in $Q_2$. (B) Competing eigenmodes explain the non-monotonic 
dynamics.}
\label{fig:coupled}
\end{figure}

\paragraph{Coupled transport.}
Two membrane transporters, $Q_1 = 
[\text{Na}^+_{\text{in}}]/[\text{Na}^+_{\text{out}}]$ and $Q_2 = 
[\text{H}^+_{\text{out}}]/[\text{H}^+_{\text{in}}]$, share the membrane 
potential with coupling matrix
\[
K = \begin{bmatrix} 1 & 0.5 \\ 0.5 & 2 \end{bmatrix}.
\]
Figure~\ref{fig:coupled}A reveals that $Q_2$ initially overshoots its steady 
state before relaxing back, while $Q_1$ monotonically approaches equilibrium. 
This non-monotonic behavior emerges from the off-diagonal coupling: the fast 
$\text{H}^+$ pump (larger $k_2$) initially races ahead, but as it depletes the 
shared energy source, it must retreat while waiting for the slower 
$\text{Na}^+$ transporter to equilibrate. Panel B reveals why: the eigenmodes 
($\lambda = 0.8, 2.2$) represent coordinated patterns where both transporters 
move together (slow mode) or in opposition (fast mode). The overshoot occurs 
when these modes interfere constructively then destructively. This example 
demonstrates that log-linear reaction quotient dynamics can generate complex 
behavior such as overshoots and non-monotonic relaxation.

\paragraph{Glycolytic oscillations.}
We demonstrate oscillatory dynamics in a simplified glycolytic system with two 
reactions: F6P $\rightleftharpoons$ FBP (phosphofructokinase) and FBP 
$\rightleftharpoons$ Products. The reaction quotients $Q_1 = 
[\text{FBP}]/[\text{F6P}]$ and $Q_2 = [\text{Products}]/[\text{FBP}]$ evolve 
according to our log-linear framework with coupling matrix:
\[
K = \begin{bmatrix} 0.5 & -2 \\ 2 & 0.5 \end{bmatrix}
\]
The off-diagonal terms capture metabolic feedback: FBP activation of 
phosphofructokinase creates positive feedback that drives oscillations.

The eigenvalues $\lambda = 0.5 \pm 2i$ have positive real part (damping rate $a 
= 0.5$ s$^{-1}$) and imaginary part (frequency $\omega = 2$ rad/s). The 
analytical solution for the reaction quotients is:
\[
Q_i(t) = K_{eq,i} \exp\left[e^{-at}(A_i\cos(\omega t) + B_i\sin(\omega 
t))\right]
\]
where $A_i$ and $B_i$ depend on initial conditions and eigenvectors. This 
yields damped oscillations with period $\tau = 2\pi/\omega \approx 3.1$ 
seconds, matching experimental observations in yeast glycolysis \cite{Dano2007}.

\begin{figure}[tb]
\centering
\includegraphics[width=\textwidth]{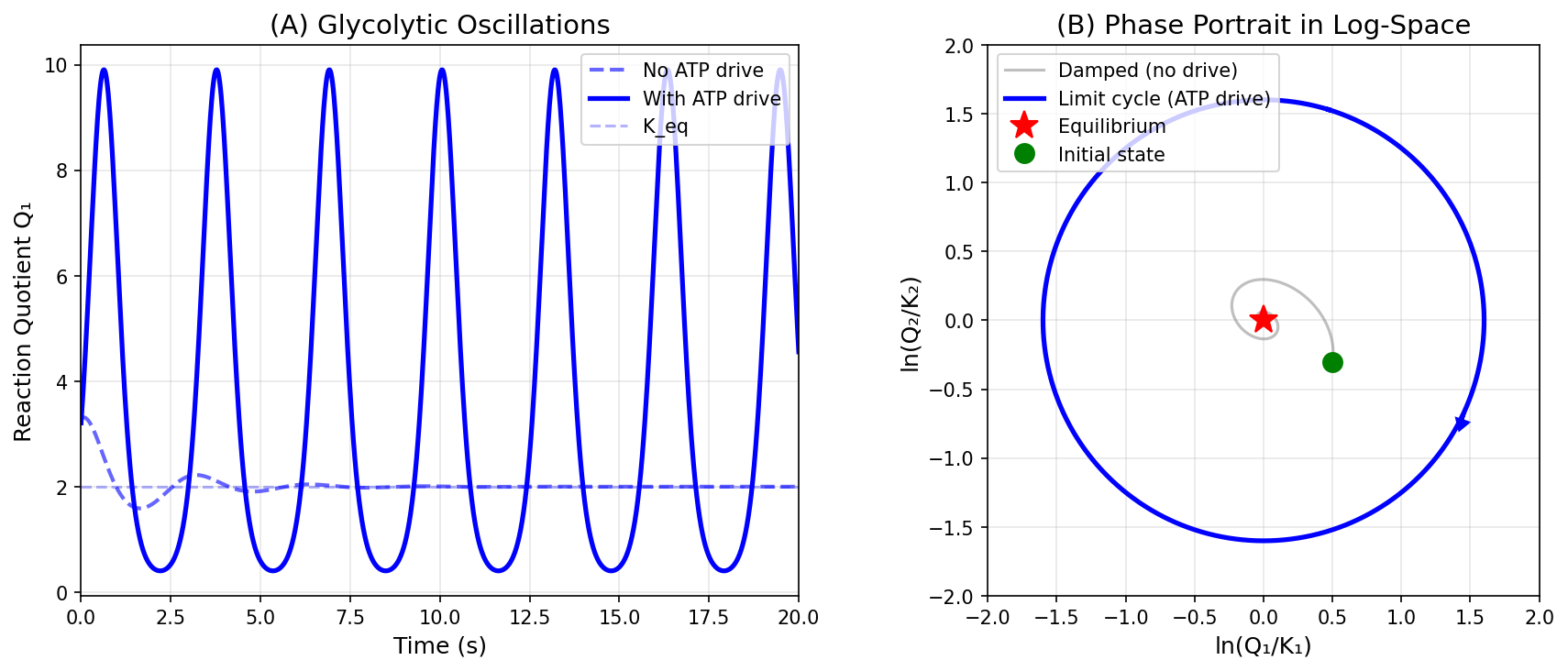}
\caption{Glycolytic oscillations in the log-linear framework. 
\textbf{(A)} Time evolution of reaction quotient $Q_1 = 
[\text{FBP}]/[\text{F6P}]$ showing damped oscillations without external drive 
(dashed line) and sustained oscillations with ATP-driven forcing at the natural 
frequency (solid line). The equilibrium value $K_{eq}$ is indicated by the gray 
dotted line. 
\textbf{(B)} Phase portrait in logarithmic coordinates $(\ln(Q_1/K_1), 
\ln(Q_2/K_2))$ demonstrating the qualitative difference between damped dynamics 
(gray spiral converging to equilibrium) and driven dynamics (blue limit cycle).}
\label{fig:glycolytic}
\end{figure}

To sustain oscillations against dissipation, we include ATP-driven forcing 
$u(t) = [u_0\sin(\omega t), 0]^T$ representing energy input from ATP 
hydrolysis. Driving at the natural frequency creates resonance, yielding a 
stable limit cycle with amplitude proportional to $u_0/a$. 
Figure~\ref{fig:glycolytic} demonstrates these dynamics: panel (A) shows how 
the reaction quotient $Q_1$ exhibits damped oscillations without external drive 
(dashed line) but maintains sustained oscillations with ATP-driven forcing 
(solid line). Panel (B) displays the corresponding phase portraits in 
log-space, where the damped system spirals to equilibrium while the driven 
system forms a stable limit cycle. Unlike mass action models requiring Hill 
functions for oscillations, our framework generates sustained metabolic rhythms 
through linear dynamics in log-space, with oscillation properties directly 
determined by the eigenvalues of $K$.

\section{Conclusion}

We have presented a framework where reaction quotients evolve exponentially 
toward equilibrium in logarithmic space, yielding linear dynamics that are 
analytically solvable. This model naturally incorporates thermodynamics through 
$\Delta G = RT\ln(Q/K_{eq})$ and decouples reaction dynamics from conservation 
constraints. Our examples demonstrate that complex behaviors such as feedback, 
oscillations, and coupled transport emerge directly from the structure of the 
coupling matrix $K$, without requiring nonlinear terms.

The analytical tractability of this framework opens several research directions.
While we have demonstrated only basic analysis here, the linear structure in 
log-space enables
the application of powerful mathematical machinery.
In metabolic engineering, control theory could optimize pathway design by 
treating $K$ as a design variable
subject to thermodynamic and biological constraints 
\cite{Stephanopoulos1998,Keasling2010}.
For drug discovery, the mathematical tractability enables systematic 
identification of therapeutic targets
and prediction of drug effects throughout the metabolic network 
\cite{Hopkins2008,Csermely2013}.
In systems medicine, eigenvalue analysis and controllability metrics from 
linear systems theory could
classify metabolic disorders and predict treatment responses 
\cite{Kitano2004,Loscalzo2011}.
The framework's linearity means that decades of control theory, from optimal 
control to robust design,
could potentially be brought to bear on problems in cellular metabolism that 
have traditionally
required complex nonlinear analysis.

While experimental validation remains essential, the log-linear framework 
offers a mathematically tractable lens for analyzing chemical reaction 
networks. Whether cells actually implement such control principles or whether 
this simply provides a useful approximation, the framework demonstrates that 
much of the complexity in chemical dynamics may be understood through linear 
algebra once viewed in the appropriate coordinates.


\clearpage

\renewcommand{\thesection}{\Alph{section}}
\appendix
\section{Conservation laws do not constrain reaction quotients}\label{app:proof}

In this appendix, we provide a self-contained proof of a well-known result from 
chemical reaction network theory (CRNT) \cite{FeinbergCRNT,HornJackson1972}: 
conservation laws impose no algebraic constraints on reaction quotients beyond 
those already determined by the stoichiometric subspace. While this 
independence is established in the CRNT literature, we include a constructive 
proof here for completeness and to make our framework accessible without 
requiring familiarity with CRNT. Our proof explicitly shows how to construct 
concentration vectors that achieve any desired reaction quotient values while 
satisfying all conservation constraints.

\paragraph{Claim.}
Let $S\in\mathbb{R}^{n\times r}$ be the stoichiometric matrix (rows: species; 
columns: reactions).
For $c\in(0,\infty)^n$, define the reaction quotients
\[
Q_j(c)=\prod_{i=1}^n c_i^{\,s_{ij}},\qquad j=1,\dots,r,
\]
stack $Q\in(0,\infty)^r$, and set $x:=\ln Q=S^\top \ln c$ (element-wise log).
Let $L\in\mathbb{R}^{n\times m}$ have columns forming a basis of $\ker S^\top$ 
and let the conserved totals be $y:=L^\top c$.
Then, for any fixed (physically attainable) totals $y^\star$, the set of 
achievable $x$ is
\[
\{x=\ln Q(c):\, c>0,\ L^\top c = y^\star\}\;=\;\operatorname{Im}S^\top .
\]

\paragraph{Proof.}
\emph{Necessity.}
For any $c>0$, write $u=\ln c$. Then
\[
x=\ln Q(c)=S^\top u\in \operatorname{Im}S^\top .
\]

\smallskip
\emph{Sufficiency.}
Fix $x_\star\in \operatorname{Im}S^\top$ and totals $y^\star$ that are 
realizable by some $c>0$ (i.e., $y^\star\in L^\top\mathbb{R}_{>0}^n$).
Choose $u_0$ with $S^\top u_0=x_\star$ and set $c_0=\exp(u_0)>0$.

\emph{Parameterization at fixed $x_\star$.}
Move multiplicatively along $\ker S^\top$:
\[
c(\alpha)\;=\;\exp(L\alpha)\circ c_0,\qquad \alpha\in\mathbb{R}^m,
\]
where $\circ$ and $\exp(\cdot)$ are element-wise. Since $S^\top L=0$,
\[
\ln Q\big(c(\alpha)\big)=S^\top\big(\ln c_0+L\alpha\big)=S^\top u_0=x_\star,
\]
so $x$ is invariant along this family.

\emph{Matching the totals.}
Define
\[
F(\alpha)\;=\;L^\top c(\alpha)\;=\;L^\top\!\big(\exp(L\alpha)\circ c_0\big),
\qquad
\Phi(\alpha)\;=\;\sum_{i=1}^n c_{0,i}\exp\!\big((L\alpha)_i\big).
\]
Then
\[
\nabla\Phi(\alpha)=F(\alpha),\qquad
\nabla^2\Phi(\alpha)=L^\top \operatorname{Diag}\!\big(\exp(L\alpha)\odot 
c_0\big)L\succ0,
\]
so $\Phi$ is $C^2$ and strictly convex (since $L$ has full column rank and the 
diagonal is positive).
Consider the strictly convex objective
\[
f(\alpha)\;=\;\Phi(\alpha)-{y^\star}^\top\alpha .
\]
Because $\Phi(\alpha)\to+\infty$ as $\|\alpha\|\to\infty$ (coercive) and $f$ is 
proper, $f$ has a unique minimizer $\alpha^\star$.
Its optimality condition is
\[
0=\nabla f(\alpha^\star)=\nabla\Phi(\alpha^\star)-y^\star
\quad\Longleftrightarrow\quad
F(\alpha^\star)=y^\star .
\]

\emph{Construction.}
Set $c^\star=\exp(L\alpha^\star)\circ c_0>0$. Then $L^\top c^\star=y^\star$ and
\[
\ln Q(c^\star)=\ln Q\big(c(\alpha^\star)\big)=x_\star .
\]

Combining necessity and sufficiency, for any feasible $y^\star$ the achievable 
$\{x=\ln Q(c):\,c>0,\ L^\top c=y^\star\}$ equals $\operatorname{Im}S^\top$, 
independent of the particular value of $y^\star$ within the attainable set. 
\hfill$\square$



\newpage
\bibliographystyle{plain} 
\bibliography{refs}       

\end{document}